\documentclass[12pt,reqno]{article}
\usepackage
[colorlinks=true, linkcolor=webgreen, filecolor=webbrown,
citecolor=webgreen]{hyperref}

\usepackage{amsmath}
\usepackage{amssymb}
\usepackage{color}
\usepackage{epsf}

\begin{document}
\makeatletter

\begin{center}
\epsfxsize=10in
\end{center}

\def\endofproofmark{$\Box$}

\begin{center}
\epsfxsize=4in
\end{center}

\def\endofproofmark{$\Box$}

\def\RR{{\mathbb R}}
\def\NN{{\mathbb N}}

\def\endofproofmark{$\Box$}

\begin{center}
\vskip 1cm {\LARGE\bf Counting and Computing by $e$} \vskip 1cm
{\large Mehdi Hassani}\\
\vskip .5cm
Department of Mathematics\\
Institute for Advanced
Studies in Basic Sciences\\
Zanjan, Iran\\
{\tt mmhassany@srttu.edu}\\
\vskip 1cm


\end{center}

\date{}

\pagestyle{myheadings}
\newtheorem{theorem}{Theorem}
\newtheorem{prop}{Proposition}
\newtheorem{lemma}{Lemma}
\newtheorem{cor}{Corollary}
\newtheorem{prob}{Note and Problem}
\def\frameqed{\framebox(5.2,6.2){}}
\def\deshqed{\dashbox{2.71}(3.5,9.00){}}
\def\ruleqed{\rule{5.25\unitlength}{9.75\unitlength}}
\def\myqed{\rule{8.00\unitlength}{12.00\unitlength}}
\def\qed{\hbox{\hskip 6pt\vrule width 7pt height11pt depth1pt\hskip 3pt}
\bigskip}
\newenvironment{proof}{\trivlist\item[\hskip\labelsep{\bf Proof}:]}{\hfill
 $\frameqed$ \endtrivlist}
\newcommand{\COM}[2]{{#1\choose#2}}

\thispagestyle{empty}
\null
\addtolength{\textheight}{2cm}

\begin{abstract}
In this paper we count the number of paths and cycles in complete
graphs by using the number $e$. Also, we compute the number of
derangements in same way. Connection by $e$ yields some nice
formulas for the number of derangements, such as
$D_n=\lfloor\frac{n!+1}{e}\rfloor$ and
$D_n=\lfloor(e+e^{-1})n!\rfloor-\lfloor en!\rfloor$, and using
these relations allow us to compute some incomplete gamma
functions and hypergeometric summations; these connections are
hidden in the heart of a nice polynomial that we call it
derangement function and a simple ordinary differential equation
concerning it.
\end{abstract}
\bigskip
\hrule
\bigskip



\bigskip
\hrule
\bigskip

\section{Introduction and Motivation}
The initial motivation of writing this paper is hidden in the
following combinatorial relations
$$
D_n=n!\sum_{i=0}^{n}\frac{(-1)^i}{i!}\hspace{5mm}(n\geq1),\\
$$
and
$$
w_n=(n-2)!\sum_{i=0}^{n-2}\frac{1}{i!}\hspace{5mm}(n\geq2),
$$
where $D_n$ is the number of derangements (permutations with no
fixed point) of $n$ distinct objects (see [3,4]), and $w_n$ is the
number of distinct paths between any pair of vertices in a
complete graph on $n$ vertices (see [4,5]). Considering
$e^{x}=\sum_{i=0}^{\infty} \frac{x^{i}}{i!}$ for $x=-1$ and $x=1$
respectively, we can get the following approximate formulas
$$
D_n\approx \frac{n!}{e} \hspace{.5cm}{\rm and}\hspace{.5cm}
w_n\approx e(n-2)!.
$$
On the other hand, it is well-known that
$$
D_n=\left\|\frac{n!}{e}\right\|\hspace{5mm}(\|x\|~\mbox{denotes
the nearest integer to }x),
$$
and this can rewritten as follows
$$
D_n=\left\lfloor\frac{n!}{e}+\frac{1}{2}\right\rfloor\hspace{5mm}(\lfloor
x\rfloor~\mbox{denotes the floor of }x).
$$
In this paper, we study these kind of combinatorial formulas
concerning partial sum of Taylor expansion of $e^x$ for some
special $x$'s. Then we use obtained results to calculate some
integrals related to the incomplete gamma function (see [1]):
$$
\Gamma(a,z)=\int_{z}^{\infty}e^{-t}t^{a-1}dt\hspace{5mm}(Re(a)>0).
$$
Also, we use them to calculate some hypergeometric summations
(see[7]):
$$
~_pF_q\left[\begin{array}{cccc}
  a_1 & a_2 & \cdots & a_p \\
  b_1 & b_2 & \cdots & b_q
\end{array};x\right]=\sum_{k\geq 0}t_k
x^k,
$$
where
$$
\frac{t_{k+1}}{t_k}=
\frac{(k+a_1)(k+a_2)\cdots(k+a_p)}{(k+b_1)(k+b_2)\cdots(k+b_q)(k+1)}x.
$$
To do these, we need some properties involving the number $e$,
that we study them in the next section.
\section{Two Interesting Formulas Involving $e$}
The obtained relations involving $e$ in this section all are
hidden in the proofs of irrationality of it. The first one is from
Rudin's analysis [8], as follows
\begin{theorem}
For every positive integer $n\geq1$, we have
\begin{eqnarray}
\sum_{i=0}^{n}\frac{n!}{i!}=\lfloor en!\rfloor.
\end{eqnarray}
\end{theorem}
\begin{proof} Since $n\geq1$, we have
$$
0<\sum_{i=0}^{\infty}\frac{1}{i!}-\sum_{i=0}^{n}\frac{1}{i!}<
\frac{1}{(n+1)!}\sum_{i=0}^\infty\frac{1}{(n+1)^i}=\frac{1}{n!n}\leq
1,
$$
and so,
$$
0<en!-\sum_{i=0}^{n} \frac{n!}{i!}\leq 1.
$$
Since $\sum_{i=0}^{n} \frac{n!}{i!}$ is an integer, the
irrationality of $e$ and the truth of the theorem both follow.
\end{proof}
The idea of the next result is hidden in Apostol's analysis [2],
where he proved the irrationality of $e$.
\begin{theorem}
Suppose $n\geq1$ is an integer, we have
$$
n!\sum_{i=0}^{n}\frac{(-1)^i}{i!}=\left\{ \begin{array}{ll}
\lfloor\frac{n!}{e}+\lambda_1\rfloor, & n ~{\rm is ~odd}, \lambda_1\in[0,\frac{1}{2}]; \\
\lfloor\frac{n!}{e}+\lambda_2\rfloor, & n ~{\rm is ~even},
\lambda_2\in[\frac{1}{3},1].
\end{array}
\right.
$$
\end{theorem}
\begin{proof} Suppose $k\geq1$ be an integer, we have
$$
0<\frac{1}{e}-\sum_{i=0}^{2k-1}\frac{(-1)^i}{i!}<\frac{1}{(2k)!},
$$
so, for every $\lambda_1$ and $\lambda_2$, we have
$$
\lambda_1<\frac{(2k-1)!}{e}+\lambda_1-\sum_{i=0}^{2k-1}\frac{(-1)^i(2k-1)!}{i!}<
\lambda_1+\frac{1}{2},
$$
and
$$
\lambda_2-1<\frac{(2k)!}{e}+\lambda_2-\sum_{i=0}^{2k}\frac{(-1)^i(2k)!}{i!}<\lambda_2.
$$
If $0\leq \lambda_1\leq\frac{1}{2}$ , then
$$
\sum_{i=0}^{2k-1}\frac{(-1)^i(2k-1)!}{i!}=\left\lfloor\frac{(2k-1)!}{e}+\lambda_1\right\rfloor.
$$
Also, if $\lambda_2\geq\frac{1}{3}$, we obtain
$$
0<\frac{(2k)!}{e}+\lambda_2-\sum_{i=0}^{2k}\frac{(-1)^i(2k)!}{i!},
$$
and with $\frac{1}{3}\leq \lambda_2\leq 1$, we yield
$$
\sum_{i=0}^{2k}\frac{(-1)^i(2k)!}{i!}=\left\lfloor\frac{(2k)!}{e}+\lambda_2\right\rfloor.
$$
This completes the proof.
\end{proof}
\begin{prob} From Theorems 1 and 2, we obtain
\begin{eqnarray}
n!\sum_{i=0}^{n}\frac{1^i}{i!}&=&\lfloor en!\rfloor\nonumber\\
n!\sum_{i=0}^{n}\frac{(-1)^i}{i!}&=&\lfloor
e^{-1}n!+\lambda\rfloor,
\hspace{5mm}\lambda\in\left[\frac{1}{3},\frac{1}{2}\right].\nonumber
\end{eqnarray}
Now, suppose $x$ is a real number, is there the set
$\Lambda_x\subseteq\mathbb{R}$, such that
$$
n!\sum_{i=0}^{n}\frac{x^i}{i!}=\lfloor
e^xn!+\lambda\rfloor,\hspace{5mm}\lambda\in\Lambda_x?
$$
for example we know $0\in\Lambda_1$ and
$[\frac{1}{3},\frac{1}{2}]\subseteq\Lambda_{-1}$.
\end{prob}
Now, we use obtained results, to study of paths and cycles in
complete graphs and then in computing the number of derangements.
\section{Paths and Cycles in Complete Graphs}
We remember that a path of length $i$ between two vertices $u$ and
$v$ is a sequence of distinct vertices as follows (see [3]),
$$
u\overset{\small(i-1)-\mbox{distindt vertices}}
{\overbrace{\bigcirc\bigcirc\bigcirc\bigcirc\cdots\bigcirc\bigcirc\bigcirc\bigcirc}}v.
$$
If $u=v$, then we have a cycle through $u$.
\begin{theorem} The number of paths between every pair of vertices in a complete graph
on $n$ vertices with $n>2$, is
$$
w_n=\lfloor e(n-2)!\rfloor,
$$
and sum of their lengths is
$$
L_w(n)=1+(n-2)\lfloor e(n-2)!\rfloor.
$$
\end{theorem}
\begin{proof} Suppose $u$ and $v$ are two vertices in a complete
graph on $n$ vertices. For counting number of paths between $u$
and $v$, we classify them according their length; the number of
paths of length $i$ is
$$
w(i)=\frac{(n-2)!}{(n-1-i)!},
$$
and since $1\leq i\leq n-1$, applying Theorem 1 with $n>2$, we
obtain
$$
w_n=\sum_{i=1}^{n-1}w(i)=\sum_{i=1}^{n-1}\frac{(n-2)!}{(n-1-i)!}=
\sum_{i=0}^{n-2}\frac{(n-2)!}{i!}=\lfloor e(n-2)!\rfloor.
$$
Also,
$$
L_w(n)=\sum_{i=1}^{n-1}iw(i)=\sum_{i=1}^{n-1}\frac{i(n-2)!}{(n-1-i)!}
=1+(n-2)\sum_{i=0}^{n-2}\frac{(n-2)!}{i!}=1+(n-2)\lfloor
e(n-2)!\rfloor.
$$
This completes the proof.
\end{proof}
\begin{cor} The average of the length of paths in a complete graph
on $n$ vertices is
$$
\frac{L_w(n)}{w_n}=n-2+O\left(\frac{1}{(n-2)!}\right),
$$
and the maximum number of paths occurs around this average.
\end{cor}
\begin{proof} According to Theorem 3, we have
$$
\frac{L_w(n)}{w_n}=n-2+\frac{1}{w_n}=n-2+O\left(\frac{1}{(n-2)!}\right).
$$
For the next assertion we have two proofs:\\
Method 1. As we mentioned in the proof of Theorem 3, the number of
paths of length $i$ is
$$
w(i)=\frac{(n-2)!}{(n-i-1)!}\hspace{5mm}(1\leq i\leq n-1),
$$
and $w(n-1)=w(n-2)=(n-2)!$; beside, for $1\leq i\leq n-3$ we have
$1\leq w(i)\leq \frac{(n-2)!}{2}$, and this yields that the
maximum of $w(i)$ occurs at $i=n+2$ and $i=n+1$; i.e, around the
average.\\
Method 2. Expand $w(i)$ as a real function with putting
$w(x)=\frac{(n-2)!}{\Gamma(n-x)}$, in which $x\in[1,n-1]$ is real.
Now, $\frac{d}{dx}w(x)=\frac{(n-2)!\Gamma'(n-x)}{\Gamma(n-x)^2}$
and so, the maximum of $w(x)$ for $1\leq x\leq n-1$, occurs
exactly at $x_{\max}=n-\gamma_0$ in which $\gamma_0\simeq
1.461632145$ is the unique zero of $\Gamma'(|x|)=0$ and therefore
$\frac{3}{2}-\gamma_0<x_{\max}-\frac{L_w(n)}{w_n}<2-\gamma_0$; so,
maximum number of paths occurs at the near of the average.
\end{proof}
The proof of the next theorem is similar to the proof Theorem 3.
\begin{theorem}
The number of cycles through every vertex in a complete graph on
$n$ vertices with $n>2$, is
$$
c_n=\lfloor e(n-1)!\rfloor-n,
$$
and sum of their lengths is
$$
L_c(n)=\lfloor en!\rfloor-\lfloor e(n-1)!\rfloor-2n+1.
$$
\end{theorem}
\begin{prob} We know that a complete graph on $n$ vertices is regular of
degree $n-1$. Can we have similar formulas about $k$-regular
graphs?
\end{prob}

\section{Computing the Number of Derangements}
Theorem 2, immediately yields the following family of formulas for
$D_n$, which is an extension of
$D_n=\lfloor\frac{n!}{e}+\frac{1}{2}\rfloor$;
\begin{theorem} For every positive integer $n\geq1$ and
$\lambda\in[\frac{1}{3},\frac{1}{2}]$, we have
$$
D_n=\left\lfloor\frac{n!}{e}+\lambda\right\rfloor.
$$
\end{theorem}
\begin{proof} We give two proofs:\\
Method 1. Using theorem 2, we have
$$
D_n=\left\{ \begin{array}{ll}
\lfloor\frac{n!}{e}+\lambda_1\rfloor, & n ~{\rm is ~odd}, \lambda_1\in[0,\frac{1}{2}]; \\
\lfloor\frac{n!}{e}+\lambda_2\rfloor, & n ~{\rm is ~even},
\lambda_2\in[\frac{1}{3},1],
\end{array} \right.
$$
and this yields the result.\\
Method 2. We know
$$
D_n=n!\left(1-\frac{1}{1!}+\cdots+\frac{(-1)^{n}}{n!}\right)
=\frac{n!}{e}+(-1)^{n}\left(\frac{1}{n+1}-\frac{1}{(n+1)(n+2)}+\cdots\right),
$$
so, for every $n\in \mathbb{N}$ we have
$$
\left|D_n-\frac{n!}{e}\right|<\frac{1}{n+1}.
$$
If $n$ is even, $D_n>\frac{n!}{e}$ and
$D_n=\lfloor\frac{n!}{e}+\lambda\rfloor$ provided
$\frac{1}{n+1}\leq \lambda\leq 1$. If $n$ is odd,
$D_n<\frac{n!}{e}$ and $D_n=\lfloor\frac{n!}{e}+\lambda\rfloor$
provided $0<\frac{1}{n+1}+\lambda\leq 1$. So we require
$\frac{1}{3}\leq\lambda\leq 1$ and $0\leq\lambda\leq\frac{1}{2}$.
This completes the second proof.
\end{proof}
\begin{cor} For every positive integer $n\geq1$, we have
\begin{eqnarray}
D_n=\left\lfloor\frac{n!+1}{e}\right\rfloor.
\end{eqnarray}
\end{cor}
\begin{proof} Use Theorem 5 with $\lambda=\frac{1}{e}$.
\end{proof}
Now, consider the idea of proving the relation
$D_n=\|\frac{n!}{e}\|$, which is the following trivial inequality
$$
\left|\frac{n!}{e}-D_n\right|\leq\frac{1}{(n+1)}+\frac{1}{(n+1)(n+2)}+\frac{1}{(n+1)(n+2)(n+3)}+\cdots,
$$
and let $M(n)$ denote the right side of above inequality. We have
$$
M(n)<\frac{1}{(n+1)}+\frac{1}{(n+1)^2}+\cdots=\frac{1}{n},
$$
and therefore,
\begin{eqnarray}
D_n=\left\lfloor\frac{n!}{e}+\frac{1}{n}\right\rfloor\hspace{5mm}(n\geq2).
\end{eqnarray}
Also, we can get a better bound for $M(n)$ as follows
$$
M(n)<\frac{1}{n+1}\left(1+\frac{1}{(n+2)}+\frac{1}{(n+2)^2}+\cdots\right)=\frac{n+2}{(n+1)^2},
$$
and similarly,
\begin{eqnarray}
D_n=\left\lfloor\frac{n!}{e}+\frac{n+2}{(n+1)^2}\right\rfloor\hspace{5mm}(n\geq2).
\end{eqnarray}
The above idea is extensible and several approximations of $M(n)$
lead us to other families of formulas for $D_n$;
\begin{theorem}
Suppose $m$ is an integer and $m\geq3$. The number of derangements
of $n$ distinct objects with $n\geq 2$ is
\begin{eqnarray}
D_n=\left\lfloor\left(\frac{\lfloor
e(n+m-2)!\rfloor}{(n+m-2)!}+\frac{n+m}{(n+m-1)(n+m-1)!}+e^{-1}\right)n!\right\rfloor
-\lfloor en!\rfloor.
\end{eqnarray}
\end{theorem}
\begin{proof}
For $m\geq3$ we have
$$
\left|\frac{n!}{e}-D_n\right|<\frac{1}{n+1}\left(1+\frac{1}{n+2}(
\cdots1+\frac{1}{n+m-1}(\frac{n+m}{n+m-1})\cdots)\right).
$$
Let $M_m(n)$ denote the right side of the above inequality; we
have
$$
M_m(n)\prod_{i=1}^{m-1}(n+i)=\frac{n+m}{n+m-1}+\sum_{j=2}^{m-1}\prod_{i=j}^{m-1}(n+i),
$$
and dividing by $\prod_{i=1}^{m-1}(n+i)$ we obtain
$$
M_m(n)=n!\left(\frac{n+m}{(n+m-1)(n+m-1)!}+\sum_{i=n+1}^{n+m-2}\frac{1}{i!}\right).
$$
Therefore,
$$
D_n=\left\lfloor\frac{n!}{e}+n!\left(\frac{n+m}{(n+m-1)(n+m-1)!}
+\sum_{i=n+1}^{n+m-2}\frac{1}{i!}\right)\right\rfloor,
$$
and reforming this by using
$\sum\limits_{i=n+1}^{n+m-2}\frac{1}{i!}=
\sum\limits_{i=0}^{n+m-2}\frac{1}{i!}-\sum\limits_{i=0}^{n}\frac{1}{i!}$,
completes the proof.
\end{proof}
\begin{cor} For $n\geq2$, we have
\begin{eqnarray}
D_n=\lfloor(e+e^{-1})n!\rfloor-\lfloor en!\rfloor.
\end{eqnarray}
\end{cor}
\begin{proof} We give two proofs:\\
Method 1. Because (5) holds for all $m\geq3$, we have
$$
D_n=\lim_{m\rightarrow\infty}\left\lfloor\left(\frac{\lfloor
e(n+m-2)!\rfloor}{(n+m-2)!}
+\frac{n+m}{(n+m-1)(n+m-1)!}+e^{-1}\right)n!\right\rfloor-\lfloor
en!\rfloor
$$
$$
=\lfloor(e+e^{-1})n!\rfloor-\lfloor en!\rfloor.
$$
Method 2. By using (1), for every $n\geq 1$ we obtain
$$
M(n)=n!(e-\sum_{i=0}^{n}\frac{1}{i!})=en!-\lfloor
en!\rfloor=\{en!\} \hspace{5mm}(\{x\}~\mbox{denotes the fractional
part of }x),
$$
and the proof follows.
\end{proof}
Now,
$$
\lim_ {m\rightarrow\infty} M_m(n)=M(n),
$$
\\
and if we put $M_1(n)=\frac{1}{n}$ and
$M_2(n)=\frac{n+2}{(n+1)^2}$ (see formulas (3) and (4)), then
$$
M_{m+1}(n)<M_m(n)\hspace{5mm}(n\geq1).
$$
Now, we find bounds sharper than $\{en!\}$ for $e^{-1}n!-D_n$ and
consequently another family of formulas for $D_n$. This family is
an extension of (6).
\begin{theorem}
Suppose $m$ is an integer and $m\geq1$. The number of derangements
of $n$ distinct objects with $n\geq2$ is
$$
D_n=\left\lfloor\left(\frac{\{e(n+2m)!\}}{(n+2m)!}+
\sum_{i=1}^{m}\frac{n+2i-1}{(n+2i)!}+e^{-1}\right)n!\right\rfloor.
$$
\end{theorem}
\begin{proof}
Since $m\geq1$, we have
\begin{eqnarray*}
\frac{e^{-1}n!-D_n}{(-1)^{n+1}}&=&
n!\sum_{i=1}^{\infty}\left(\frac{1}{(n+2i-1)!}-\frac{1}{(n+2i)!}\right)\\&<&
n!\left(\sum_{i=1}^{m}\frac{n+2i-1}{(n+2i)!}+\sum_{i=2m+1}^{\infty}\frac{1}{(n+i)!}\right).
\end{eqnarray*}
Let $N_m(n)$ denote the right member of above inequality.
Considering (1), we obtain
$$
N_m(n)=n!
\sum_{i=1}^{m}\left(\frac{n+2i-1}{(n+2i)!}+\frac{\{e(n+2m)!\}}{(n+2m)!}\right),
$$
and for $n\geq2$, we yield that $D_n=\lfloor
e^{-1}n!+N_m(n)\rfloor$. This completes the proof.
\end{proof}
\begin{cor} For all integers $m, n\geq1$, we have
$$
N_{m+1}(n)<N_m(n),\ \ \ N_1(n)<\{en!\}.
$$
\end{cor}
Therefore we have the following chain of bounds for
$|\frac{n!}{e}-D_n|$
$$
\left|\frac{n!}{e}-D_n\right|<\cdots<N_2(n)<N_1(n)<\{en!\}<\cdots<M_2(n)<M_1(n)<1\hspace{5mm}(n\geq2).
$$
In the next section you will see that how searching other formulas
for $D_n$ lead us to connections between $D_n$, incomplete gamma
functions and hypergeometric summations.

\section{Applications: Derangement Function, Incomplete
Gamma Function and Hypergeometric Function}

Let's find other formulas for $D_n$. The computer algebra program
MAPLE yields that
$$
D_n=(-1)^n\texttt{hypergeom}([1,-n],[~~],1),
$$
and
$$
D_n=e^{-1}\Gamma(n+1,-1),
$$
where $\texttt{hypergeom}([1,-n],[~~],1)$ is MAPLE's notation for
a hypergeometric function; more generally,
$$
\texttt{hypergeom}([a_1~~a_2~\cdots~a_p],[b_1~~b_2~\cdots~b_q],x)=
~_pF_q\left[\begin{array}{cccc}
  a_1 & a_2 & \cdots & a_p \\
  b_1 & b_2 & \cdots & b_q
\end{array};x\right].
$$
Now, because we know the value of $D_n$, we can estimate some
summations and integrals. To do this, we define the
\textit{derangement function}, a natural generalization of the
number of derangements, denoted by $D_n(x)$, for every integer
$n\geq0$ and every real $x$ as follows
$$
D_n(x)=\left\{ \begin{array}{ll}
n!\sum_{i=0}^{n}\frac{x^i}{i!}, & x\neq0; \\
n!, & x=0.
\end{array}
\right.
$$
It is easy to obtain the following generalized recursive relations
$$
D_n(x)=(x+n)D_{n-1}(x)-x(n-1)D_{n-2}(x)=x^n+nD_{n-1}(x),\hspace{5mm}
(D_0(x)=1,D_1(x)=x+1).
$$
The function $D_n(x)$ connects proven results in two previous
sections with incomplete gamma functions and hyperbolic
summations. Also, these connection has some consequences, which we
collected some of them in the next theorem and corollaries.
\begin{theorem} Suppose $x$ is a real number, we have
\begin{eqnarray}
D_n(x)=x^n~_2F_0\left[\begin{array}{cc}
  1 & -n \\
  -
\end{array};-\frac{1}{x}\right]\hspace{5mm}(x\neq 0),
\end{eqnarray}
and
\begin{eqnarray}
D_n(x)=e^x\Gamma(n+1,x).
\end{eqnarray}
\end{theorem}
\begin{proof} The function $D_n(x)$ satisfies in the following
differential equation
$$
D_n(x)-\frac{d}{dx}D_n(x)=x^n,\hspace{5mm}D_n(0)=n!,
$$
solving this equation by summations yields (7) and solving it by
integrals yields (8).
\end{proof}
\begin{cor} For every $n\in\mathbb{N}$, we have
$$
~_2F_0\left[\begin{array}{cc}
  1 & -n \\
  -
\end{array};-1\right]=\lfloor en!\rfloor,
$$
and
$$
~_2F_0\left[\begin{array}{cc}
  1 & -n \\
  -
\end{array};1\right]=
(-1)^n\left\lfloor\frac{n!+1}{e}\right\rfloor.
$$
\end{cor}
\begin{proof} Consider (7) with $x=-1$ and $x=1$ respectively, and
use the relations (1) and (2).
\end{proof}
\begin{cor}
For every real $x\neq0$, we have
$$
~_1F_1\left[\begin{array}{c}
  n+1 \\
  n+2
\end{array};-x\right]=\frac{(n+1)(n!-e^{-x}D_n(x))}{x^{n+1}}.
$$
\end{cor}
\begin{proof}
Obvious.
\end{proof}
\begin{cor}
For every integer $n\geq1$ we have
$$
\int_{-1}^{\infty}e^{-t}t^ndt=e\left\lfloor\frac{n!+1}{e}\right\rfloor,
$$
$$
\int_{0}^{\infty}e^{-t}t^ndt=n!,
$$
$$
\int_{1}^{\infty}e^{-t}t^ndt =\frac{\lfloor en!\rfloor}{e},
$$
and
$$
\int_{0}^{1}e^{-t}t^ndt =\frac{\{en!\}}{e},
$$
$$ \int_{-1}^{0}e^{-t}t^ndt =\left\{
\begin{array}{cc}
-e\{\frac{n!}{e}\} & n ~{\rm is ~odd},\\
e-e\{\frac{n!}{e}\} & n ~{\rm is ~even},
\end{array}
\right.
$$
$$
\int_{-1}^{1}e^{-t}t^ndt=e\lfloor(e+e^{-1})n!\rfloor-(e+e^{-1})\lfloor
en!\rfloor,
$$
\end{cor}
\begin{proof} Consider the definition of incomplete gamma function
and use relation (8) with $x=-1$, $x=0$ and $x=1$ respectively,
and use (1), (5) and (6).
\end{proof}

\begin{prob} Note that $D_n(x)$ is a nice polynomial.
Its value for $x=-1$ is $D_n$, for $x=0$ is the number of
permutations of $n$ distinct objects and for $x=1$ is $w_{n+2}=$
the number of distinct paths between every pair of vertices in a
complete graph on $n+2$ vertices. Is there any combinatorial
meaning for the value of $D_n(x)$ for other values of $x$?
\end{prob}

\begin{prob} Since, $n!=\Gamma(n+1)$ in which
$\Gamma(x)=\int_{0}^{\infty}e^{-t}t^{x-1}dt$, we can generalize
the derangement function as follows
$$
D_\lambda(x)=\left\{ \begin{array}{ll}
\int_{0}^{\lambda}\frac{x^tdt}{\Gamma(t+1)}, & x\neq0; \\
\Gamma(\lambda-1), & x=0,
\end{array}
\right.
$$
in which $\lambda>0$ is a real number. Here our famous $e$ can be
replace by $E=\int_{0}^{\infty}\frac{dt}{\Gamma(t+1)}\simeq
2.266534508$. Can we study $D_\lambda(x)$ and yield some new nice
results?
\end{prob}

\end{document}